\documentclass[12pt]{article}
\newcommand{\copyleft}{
GNU FDL\thanks{
Copyright (C) 2006 J. Laurie Snell.
Permission is granted to copy, distribute and/or modify this document
under the terms of the GNU Free Documentation License, 
as published by the Free Software Foundation;
with no Invariant Sections, no Front-Cover Texts, and no Back-Cover Texts.
}}
\title{The Engel algorithm for absorbing Markov chains}
\author{
J. Laurie Snell
\thanks{
This module was produced under a project directed
by J.\ Laurie Snell and funded by the National Science Foundation.
The purpose of this project was to provide modules to
serve as supplementary material for courses in probability and,
where appropriate, to illustrate the use of the computer in probability.
This module was resurrected through the efforts of Peter Doyle.}
}
\date{Version dated circa 1979
\\ \copyleft
}

\usepackage{graphics, epsfig}
\newcommand{\putfig}[3]{
\begin{figure}[cbt]
\centerline{\mbox{\includegraphics*[width=5.0truein]{figures/#2.eps}}}
\caption{#3}
\label{fig:#1}
\end{figure}
}
\newcommand{\figref}[1]{\ref{fig:#1}}

\newdimen\snellbaselineskip
\newdimen\snellskip
\snellbaselineskip=\baselineskip
\def\srule{\omit\kern.5em\vrule\kern-.5em}
\newbox\bigstrutbox
\setbox\bigstrutbox=\hbox{\vrule height14.5pt depth9.5pt width0pt}
\def\bigstrut{\relax\ifmmode\copy\bigstrutbox\else\unhcopy\bigstrutbox\fi}
\def\middlehrule#1#2{\noalign{\kern-\snellbaselineskip\kern\snellskip}
&\multispan#1\strut\hrulefill
&\omit\hbox to.5em{\hrulefill}\vrule
height \snellskip\kern-.5em&\multispan#2\hrulefill\cr}

\makeatletter
\def\bordermatrix#1{\begingroup \m@th
  \@tempdima 8.75\p@
  \setbox\z@\vbox{%
    \def\cr{\crcr\noalign{\kern2\p@\global\let\cr\endline}}%
    \ialign{$##$\hfil\kern2\p@\kern\@tempdima&\thinspace\hfil$##$\hfil
      &&\quad\hfil$##$\hfil\crcr
      \omit\strut\hfil\crcr\noalign{\kern-\snellbaselineskip}%
      #1\crcr\omit\strut\cr}}%
  \setbox\tw@\vbox{\unvcopy\z@\global\setbox\@ne\lastbox}%
  \setbox\tw@\hbox{\unhbox\@ne\unskip\global\setbox\@ne\lastbox}%
  \setbox\tw@\hbox{$\kern\wd\@ne\kern-\@tempdima\left(\kern-\wd\@ne
    \global\setbox\@ne\vbox{\box\@ne\kern2\p@}%
    \vcenter{\kern-\ht\@ne\unvbox\z@\kern-\snellbaselineskip}\,\right)$}%
  \null\;\vbox{\kern\ht\@ne\box\tw@}\endgroup}

\makeatother

\begin{document}
\maketitle
\begin{abstract}
In this module,
suitable for use in an introductory probability course,
we present Engel's chip-moving algorithm
for finding the basic descriptive quantities
for an absorbing Markov chain,
and prove that it works.
The tricky part of the proof involves showing that the initial distribution
of chips recurs.
At the time of writing (circa 1979)
no published proof of this was available,
though Engel had stated that such a proof had been found
by L.\ Scheller.

\vspace{0.5cm}
\end{abstract}

In \cite{engel:text} Arthur Engel provides an interesting algorithm
for finding the basic descriptive quantities for an absorbing Markov chain.
It was developed as a teaching tool.
It can be carried out by moving chips around the transition graph
of the Markov chain.
In \cite{engel:why} Engel provides examples and a partial proof that the
algorithm works.
We complete his proof here.

Let $P$ be the transition matrix of an absorbing Markov chain
with states $1,2,\ldots,r$ transient and states $r+1,\ldots,s$
absorbing.
Assume that the transition probabilities are rational numbers
and write them in the form $p_{ij} = \frac{r_{ij}}{r_i}$
where $r_i$ and $r_{ij}$ are integers.
Note that if we have $r_i$ chips in state $i$,
we can distribute these chips according to $P$ by sending $r_{ij}$
to state $j$.

We choose a particular transient state $u$.
The algorithm will find $N_{uj}$, the expected number of visits to
state $j$ when the process is started in state $u$,
and $B_{uj}$, the probability that the process ends in
absorbing state $j$ when it is started in state $u$.
We put an initial disbribution of chips
$c_1,c_2,\ldots,c_r$ on the transient states.
The choice $c_i = r_i-1$ for $i \neq u$
and $c_u=r_u$ is called a \emph{critical loading}.
We put no chips on the absorbing states initially.
We are allowed to make two types of moves.
A \emph{type 1 move} is to move $r_i$ chips from
state $i$ by moving $r_{ij}$ to state $j$ for each $j$.
A \emph{type 2 move} is to add a chip to our starting state $u$.
A type 1 move is possible only when there are at least $r_i$ chips
on state $i$.
We also assume that we are only allowed to make a type 2 move
when no type 1 move is possible.
We begin with a critical loading and then make as many type 1 moves
as possible.
When no more such moves are possible, we employ type 2 moves until
we have $r_u$ chips at $u$.
We then see if we have again a critical loading.
If not, we continue the process.
We stop when we obtain a critical loading again.
We shall prove that this will occur.

Let $w_{uj}$ be the total number of chips that have been moved
out of the transient state $j$ and $v_{uk}$ the number of chips
that have been moved into the absorbing state $k$ during the
algorithm.
Let $v_u = \sum_k v_{uk}$.
Then we shall show that
(a) the algorithm will stop, and
(b) $N_{uj} = \frac{w_{uj}}{v_u}$
and $B_{uj} = \frac{v_{uj}}{v_u}$.

First we shall illustrate the algorithm by a simple example.
We consider random walk with a drift with state 0 and 3 absorbing.
(See Figure \figref{chain}.)
\putfig{chain}{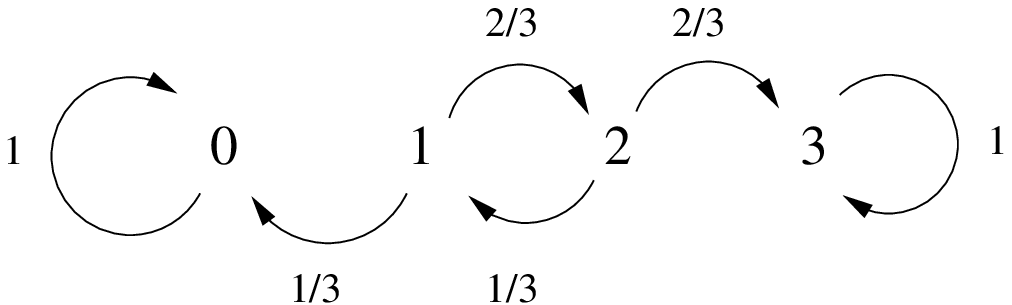}{Random walk with a drift.}
We choose $u$ to be state 1.
The critical loading is $(3,2)$.
Carrying out the algorithm (see Table \ref{engelrun}),
\begin{table}
\caption{A run of Engel's algorithm}
\label{engelrun}
\[
\begin{array}{ccccc|c}
\mbox{State}&0&1&2&3&\mbox{Move}\\
\hline
\mbox{Start}&0&3&2&0&1_1\\
&1&0&4&0&1_2\\
&1&1&1&2&2\\
&1&2&1&2&2\\
&1&3&1&2&1_1\\
&2&0&3&2&1_2\\
&2&1&0&4&2\\
&2&2&0&4&2\\
&2&3&0&4&1_1\\
&3&0&2&4&2\\
&3&1&2&4&2\\
&3&2&2&4&2\\
\mbox{Final}&3&3&2&4\\
\end{array}
\]
\end{table}
we obtain
$w_{11}=9$,
$w_{12}=6$,
$v_{10}=3$,
$v_{13}=4$, and
$v_1=7$.
Thus
$N_{11}=9/7$,
$N_{12}=6/7$,
$B_{10}=3/7$, and
$B_{13}=4/7$.

Let us carry out these calculations by the standard method.
We put $P$ in canonical form

\[
\offinterlineskip
P\;= \bordermatrix{
                               &\hbox{TR.}&\omit\hfil&\hbox{ABS.}\cr
           \hbox{TR.}\bigstrut &Q   &\srule    &R    \cr
\middlehrule{1}{1}
           \hbox{ABS.}\bigstrut&0   &\srule    &I}
:
\]
\[
\offinterlineskip
P\;= \bordermatrix{&
               \hbox{1} &\hbox{2}&\omit\hfil&\hbox{0}&\hbox{3}\cr
\hbox{1}\strut  &  0    &2/3 &  \srule  & 1/3    &  0     \cr
\hbox{2}\strut  &1/3    &  0 &  \srule  & 0      &  2/3   \cr
\middlehrule{2}{2}
\hbox{0}\strut  &  0    &  0 &  \srule  & 1      &  0     \cr
\hbox{3}\strut  &  0    &  0 &  \srule  & 0      &  1}\
.
\]
Then
\[
N = (I-Q)^{-1} =
\left(\begin{array}{cc}1&-2/3\\-1/3&1\end{array}\right)^{-1}
=
\bordermatrix{&1&2\cr 1&9/7&6/7\cr 2&3/7&9/7}
\]
and
\[
B = NR =
\left(\begin{array}{cc}9/7&6/7\\3/7&9/7\end{array}\right)
\left(\begin{array}{cc}1/3&0\\0&2/3\end{array}\right)
=
\bordermatrix{&0&3\cr 1&3/7&4/7\cr 2&1/7&6/7}
.
\]

We turn now to the proofs.
Let us start with any initial loading less than or equal to the critical
loading.
We continue the algorithm until some distribution
$b_1,b_2,\ldots,b_r$ 
is repeated.
There must be such a distribution since there are only a finite
number of distributions that can occur.

We now start over again, this time with the initial loading
$b_1,b_2,\ldots,b_r$.
We shall follow the proof in \cite{engel:why}
to show that our foumulas for $N _{uj}$ and $B_{uj}$ are correct.
We shall use only the fact that the initial distribution repeats.
We shall then prove that the critical loading will always repeat.

Since the initial and final distributions are the same,
the number of chips that move into a state must equal the number
that move out.
Thus
\[w_{uu} = v_u + \sum_k w_{uk}P_{ku}
\]
\[
w_{ui} = \sum_k w_{uk}P_{ki} \;\;\mbox{for $i \neq u$}
.
\]
Let $\bar{w}_{uj} = \frac{w_{uj}}{v_u}$.
Then we have
\[
\bar{w}_{uu} = 1 + \sum_k \bar{w}_{uk}P_{ku}
\]
\[
\bar{w}_{ui} = \sum_k \bar{w}_{uk}P_{ki} \;\;\mbox{for $i \neq u$}
.
\]
If we do this for all choices of $u$,
we obtain the matrix equations
\[
\bar{W} = I + \bar{W}Q
\]
or
\[
\bar{W}(I-Q) = I
\]
and
\[
\bar{W} = (I-Q)^{-1} = N
.
\]
Similarly for $j$ absorbing
\[
v_{uj}=\sum_{k} w_{uk}P_{kj}
.
\]
Putting $\bar{v}_{uj} = \frac{v_{uj}}{v_u}$,
we obtain
\[
\bar{v}_{uj} = \sum_k \bar{w}_{uk}P_{kj}
.
\]
Doing this for all $u$, we obtain in matrix form
\[
\bar{V} = \bar{W}R
.
\]
But since $\bar{W}=N$, we have
\[
\bar{V} = NR = B
.
\]

We show now that the critical loading will repeat.
In \cite{engel:why} it is stated that L.\ Scheller has proved this fact
but the proof is not given.
The following proof was provided by Peter Doyle.

We know that there is an initial sequence
$b_1,b_2,\ldots,b_r$
that repeats.
We can assume that $b_i<r_i$ for $i \neq u$ and $b_u = r_u$.
We choose this sequence using green chips.
We then add red chips to bring this distribution up to the critical loading.
We now carry out the same sequence of moves using only the green chips.
Then when our initial distribution repeats,
we will have the critical loading again since we have not touched the red
chips.
Let us call this modified procedure \emph{method A}.
We call it `modified' because we have not followed the Engel rules
for our initial load---the critical loading.
Assume that in method A we added $m$ chips to state $u$.
Then we consider a second method which we call \emph{method B}.
In method B, we start with the critical loading
and carry out the Engel rules until we have added $m$ chips
and can make no more type 1 moves.
We will show that the two methods lead to the same final distribution.
Since we know that the final distribution in method A is the critical loading,
the same will be true for method B.

Let
$e_1,e_2,\ldots,e_n$
be the moves in order as made in method A and
$f_1,f_2,\ldots,f_\nu$
those made in method $B$.
Moves $e_1$ and $f_1$
are the same,
namely move $r_u$ chips from $u$.
We shall prove by induction that all the $e$ moves appear somewhere
in the $f$ list.
Assume that we have, for the moves
$e_1,e_2,\ldots,e_k$,
found corresponding moves
$f_{i_1},f_{i_2},\ldots,f_{i_k}$ among the $f$'s.
Consider move $e_{k+1}$.
If this is a type 2 move,
there must be a corresponding move $f_{i_{k+1}}$
since we make the same number
$m$ of such moves in these two methods.
Assume then that $e_{k+1}$ is a type 1 move,
say move $r_t$ chips from state $t$.
If there is such a move $f_{i_{k+1}}$ among the $f$'s
which occurs in addition to
$f_{i_1},f_{i_2},\ldots,f_{i_k}$
and before these are all carried out,
we can choose this move to correspond to $e_{k+1}$.
Assume that there is no such move.
The move $e_{k+1}$ was made possible by the moves 
$e_1,e_2,\ldots,e_k$.
This move must also be possible by virtue of the moves
$f_{i_1},f_{i_2},\ldots,f_{i_k}$.
The other moves that could have occurred while these moves
are made can only help by adding chips to state $t$.
Thus again, we can match $e_{k+1}$ and the induction step is proved.

The same argument show that all the $f_k$'s are among the $e_k$'s.
Thus the two sets of moves are the same,
and this implies that methods $A$ and $B$ lead to the same final
distribution, as we wanted to show.

\bibliography{engel}
\bibliographystyle{plain}

\end{document}